\documentclass[a4paper,12pt]{article}
\usepackage[english]{babel}
\usepackage{amsfonts,amssymb,amsmath}

\begin{document} 
 
 \centerline{\bf\Large Two-parameter asymptotics in the Cauchy problem}
 \vskip 1mm  
 \centerline{\bf\Large for a parabolic equation} 

 \vskip 5mm 
 \centerline{\bf S.V.~Zakharov} 

 \vskip 5mm 

\begin{center}
Institute of Mathematics and Mechanics,\\
Ural Branch of the Russian Academy of Sciences,\\
16, S.Kovalevskaja street, 620990, Ekaterinburg, Russia
\end{center}

\vskip 5mm 

\textbf{Abstract.}
The Cauchy problem for a quasi-linear
parabolic equation with a small parameter
at a higher derivative is considered. 
The initial step-like function
contains another small parameter. 
Formal asymptotic solutions of the problem
in small parameters are constructed.

\vskip 3mm 

 
\vskip 5mm

 \section{Introduction} 
 
 In the present work, we consider the Cauchy problem
for a quasi-linear parabolic equation: 
 \begin{eqnarray} 
 \label{eq} 
 \frac{\partial u}{\partial t} + 
 \frac{\partial \varphi(u)}{\partial
x} = 
 \varepsilon \frac{\partial^2 u}{\partial
x^2},& 
 \quad t\geqslant 0, & \quad \varepsilon>0, \\ 
 \label{ic}\phantom{\frac{1}{1}} 
 u(x,0,\varepsilon,\rho) = \nu ( {x}{\rho}^{-1}),& 
 \quad x\in\mathbb{R}, & \quad \rho>0. 
 \end{eqnarray} 
 We assume that $\varepsilon>0,$ 
the function $\varphi$ is infinitely differentiable
 and its second derivative is strictly positive. 
 The initial function~$\nu$ is bounded and smooth. 
 
This model is used for studying 
the evolution of a wide class of physical systems 
with a small dissipation and probabilistic
processes~\cite{bu,wh,ib}. 
The interest to problem under consideration is explained
by physical applications and 
 the fact that its solutions allow one to obtain viscous 
generalized solutions of  the limit equation. 
This problem had been studied by
N.S.~Bakhvalov, I.M.~Gelfand, A.M.~Il'in, 
E.~Hopf, O.A.~Ladyzhenskaya, O.A.~Oleinik,
 and many other mathematicians. 
It is strictly proved~\cite{lsu}, 
 that there exists a unique bounded infinitely
differentiable with respect to $x$ and~$t$
solution $u(x,t,\varepsilon)$. 
 
 The aim of the present paper is to construct 
asymptotics solutions $u(x,t,\varepsilon,\rho)$ 
of problem~(\ref{eq})--(\ref{ic}) 
 as $\varepsilon\to 0$ and $\rho\to 0$. 
The structure of asymptotic series
essentially depends
on the relation between parameters $\varepsilon$ and $\rho$
as shown below. 
 
 First, note that the change of variables $$ 
 x = \rho \sigma, \qquad t = \rho \theta 
 $$ 
 in equation~(\ref{eq}) and the initial
condition~(\ref{ic}) 
 leads to the following problem: 
 $$ 
 \frac{\partial u}{\partial \theta} + 
 \frac{\partial \varphi(u)}{\partial \sigma} = 
 \delta \frac{\partial^2 u}{\partial \sigma^2}, 
 \qquad 
 \delta = \frac{\varepsilon}{\rho}, 
 \qquad 
 u(\sigma,0) = \nu (\sigma). 
 $$ 
 An asymptotic approximation
 of the solution of such a problem
up to an arbitrary power of parameter $\delta$ 
 is obtained directly from~\cite[Ch. VI]{ib}. 
 For example,
 the expansion of the solution in a neighborhood
of the singular point $(x,t)=(0,\rho)$ 
 has the form $$ 
 \sum\limits_{k=1}^{\infty} \delta^{k/4} 
 \sum\limits_{j=0}^{k-1} w_{k,j}(\xi,\tau) \ln^j \delta, 
 $$ 
 where coefficients $w_{k,j}(\xi,\tau)$ depend on the
inner variables, which are determined
using change $$ 
 x = \varepsilon^{3/4} \rho^{1/4} \xi, 
 \qquad 
 t = \rho + \varepsilon^{1/2} \rho^{1/2} \tau. 
 $$ 
 The leading term of this expansion
is found with the help of
the Cole--Hopf transform:
$$ 
 w_{1,0}(\xi,\tau)= 
 - \frac{2}{\varphi''(0) \Lambda(\xi,\tau)} 
 \frac{\partial \Lambda(\xi,\tau)}{\partial \xi},
 $$ 
 $$ 
 \Lambda(\xi,\tau) = \int\limits_{-\infty}^{\infty} 
 \exp ( -2 z^4 + z^2 \tau + z \xi) \,dz. 
 $$ 
 
 Now, let the relation between parameters ${\varepsilon}$
and ${\rho}$ be  such that
 $$\displaystyle\frac{\rho}{\varepsilon}\to 0.$$ 
Assume that the function $\nu$ in~(\ref{ic})
 satisfies the following asymptotic relations: 
 \begin{equation}\label{na} 
 \nu(\sigma) = \sum\limits_{n=0}^{\infty} 
 \frac{\nu^{\pm}_n}{\sigma^{n}}, 
 \qquad 
 \sigma\to \pm\infty. 
 \end{equation} 

 In papers~\cite{2ps,zz}, it is shown
 that for the solution
of problem $(\ref{eq})$--$(\ref{ic})$ 
 as $\varepsilon\to 0$, $\rho\to 0$, 
 $\mu=\rho/\varepsilon\to 0$ 
 in the strip $$ 
 \{ (x,t) : x\in\mathbb{R},\ 0 \leqslant
t\leqslant T\} 
 $$ 
 there holds the asymptotic formula
 $$ 
 {u}(x,t,\varepsilon,\rho)= 
 h_0\left( 
 \frac{x}{\rho}, 
 \frac{\varepsilon t}{\rho^2}\right)- R_{0,0,0}\left( 
 \frac{x}{2\sqrt{\varepsilon t}} 
 \right)+ 
 \Gamma\left( 
 \frac{x}{\varepsilon}, 
 \frac{t}{\varepsilon} 
 \right) 
 +O(\mu^{1/2}\ln\mu), 
 $$ 
 where $$ 
 h_0(\sigma,\omega) = 
 \frac{1}{2\sqrt{\pi\omega}} 
 \int\limits_{-\infty}^{+\infty} 
 \nu(s) \exp\left[ -\frac{(\sigma-s)^2}{4\omega}\right]ds, 
 $$ 
 $$ 
 R_{0,0,0}(z) = \frac{\nu^{-}_{0}}{\sqrt{\pi}} 
 \int\limits_{z}^{+\infty} \exp(-y^2)\,
dy 
 + \frac{\nu^{+}_{0}}{\sqrt{\pi}} 
 \int\limits_{-\infty}^{z} \exp(-y^2)\,
dy, 
 $$ 
 $\Gamma$ is the solution of the equation $$ 
 \frac{\partial \Gamma}{\partial \theta} + 
 \frac{\partial \varphi(\Gamma)}{\partial \eta} - \frac{\partial^2 \Gamma}{\partial \eta^2} = 0 
 $$ 
 in the inner variables $\eta = x/\varepsilon$, 
 $\theta = t/\varepsilon$ 
 with the initial condition
 $$ 
 \Gamma(\eta,0) = 
 \begin{cases} 
 \nu^-_0, & \eta<0, \\ 
 \nu^+_0, & \eta\geqslant 0. 
 \end{cases} 
 $$ 
 
In the present paper, formal asymptotic
solutions of problem $(\ref{eq})$--$(\ref{ic})$
are constructed in the form of infinite series.

 \setcounter{equation}{0} 
 \section{The outer expansion} 
 
 The behavior of the solution 
of problem $(\ref{eq})$--$(\ref{ic})$ is mainly 
determined by the solution of the limit problem 
\begin{equation}\label{lcp} 
 \frac{\partial u}{\partial t} + 
 \frac{\partial \varphi(u)}{\partial
x} =0, 
 \qquad 
 u(x,0) = 
 \begin{cases} 
 \nu^-_0, & x<0, \\ 
 \nu^+_0, & x\geqslant 0. 
 \end{cases} 
 \end{equation} 
 For $\nu^-_0>\nu^+_0$ using the method 
of characteristics,
 we find its generalized solution $$ 
 u_{0,0}(x,t) = 
 \begin{cases} 
 \nu^-_0, & x<ct, \\ 
 \nu^+_0, & x>ct, 
 \end{cases} 
 \qquad 
 {c}= 
 \frac{\varphi(\nu_0^+)-\varphi(\nu_0^-)}{\nu_0^+ -\nu_0^-}. 
 $$ 
This solution is discontinuous on the line of the shock wave $x=ct$.

 First, let us find the outer expansion
 in the domain
 $$ 
 \Omega^+_0 = \{ (x,t)\,: \, x>ct+\varepsilon^{1-\delta_0}, 
 \ 0<\delta_0<1\}.
 $$ 
 Taking into account~(\ref{na}), we will construct
the outer asymptotic expansion in the form of the series 
\begin{equation}\label{outp} 
 U^+(x,t,\varepsilon,\rho) = \nu^+_0 + 
 \sum\limits_{m = 1}^{\infty}\sum\limits_{n=0}^{m-1} 
 \rho^{m-n}\varepsilon^{n} u^+_{m,n}(x,t). 
 \end{equation} 
 In the domain $$ 
 \Omega^-_0 = \{ (x,t)\,: \, x<ct-\varepsilon^{1-\delta_0}\} 
 $$ 
 we will construct an analogous series 
\begin{equation}\label{outm} 
 U^-(x,t,\varepsilon,\rho) = \nu^-_0 
 + \sum\limits_{m = 1}^{\infty}\sum\limits_{n=0}^{m-1} 
 \rho^{m-n}\varepsilon^{n} u^-_{m,n}(x,t). 
 \end{equation} 
 Formally substituting series~(\ref{outp})
and~(\ref{outm}) 
 into equation~(\ref{eq}) and
collecting coefficients at $\rho^{m-n}\varepsilon^{n}$, 
 we arrive at the recurrence system of initial value
problems 
\begin{equation}\label{os} 
 \frac{\partial u^{\pm}_{m,n}}{\partial
t} 
 + \varphi'(\nu^{\pm}_0) 
 \frac{\partial u^{\pm}_{m,n}}{\partial
x} 
 =F^{\pm}_{m,n}, 
 \quad 
 u^{\pm}_{m,n}(x,0)= \delta_{n,0}\nu^{\pm}_{m}
x^{-m}, 
 \end{equation} 
 where $\delta_{0,0}=1$, $\delta_{n,0}=0$
for $n\neq 0$, 
 \begin{equation}\label{fmnp} 
 F^{\pm}_{m,n} = 
 \frac{\partial^2 u^{\pm}_{m-1,n-1} }{\partial
x^2} 
 - \sum\limits_{q=2}^{m-n} 
 \frac{\varphi^{(q)}(\nu^{\pm}_0)}{q!} 
 \sum\limits_{\substack{i_1+\dotsc+i_q = m\\
j_1 +\dotsc+ j_q = n}} 
 \frac{\partial \left(u_{i_1,j_1}\cdot \dotsc\cdot
u_{i_q,j_q} \right) } 
 {\partial x}. 
 \end{equation} 
Using the method of characteristics,
 we find the coefficients of the outer expansion: 
 \begin{equation}\label{upmn} 
 u^{\pm}_{m,n}(x,t) = 
 \frac{\delta_{n,0}\nu^{\pm}_{m} }{[x-\varphi'(\nu^{\pm}_0)t]^{m}} 
 + \int\limits_{0}^{t} 
 F^{\pm}_{m,n}(x-\varphi'(\nu^{\pm}_0)(t-t'),t')
dt'. 
 \end{equation} 
 Thus, for $t=0$ formal series~(\ref{outp})
and~(\ref{outm}) 
become asymptotic series for the initial function: 
 $$ 
 U^{\pm}(x,0,\varepsilon,\rho)=\sum\limits_{m=0}^{\infty} 
 \nu^\pm_m \left(\frac{\rho}{x} \right)^m. 
 $$ 
 
 Using relations~(\ref{fmnp}) and~(\ref{upmn}),
 by induction we arrive at the following statement. 

\textbf{Theorem~1.}
{\it For $m\geqslant 1$ and $0\leqslant n\leqslant m-1$,
there holds formula 
\begin{equation}\label{umn} 
 u^{\pm}_{m,n}(x,t) = 
 \sum\limits_{s=n}^{m-1} 
 \frac{ \alpha^{\pm}_{m,n,s}\, t^s}{[x-\varphi'(\nu^{\pm}_0)t]^{m+s}}, 
 \end{equation} 
 {where} $\alpha^{\pm}_{m,n,s}$ are constants. 
}

 \setcounter{equation}{0} 
 \section{The inner expansion} 

We make the change of variables 
\begin{equation}\label{cv} 
 x = \rho\sigma, 
 \qquad 
 t = \frac{\rho^2}{\varepsilon} \omega. 
 \end{equation} 
 Then equation~(\ref{eq}) becomes 
\begin{equation}\label{ea} 
 \frac{\partial h}{\partial \omega} - \frac{\partial^2
h}{\partial \sigma^2} = 
 - \mu \frac{\partial \varphi(h)}{\partial \sigma}, 
 \end{equation} 
 where $h(\sigma,\omega)\equiv u(\rho\sigma,\rho^2 \omega/\varepsilon)$, 
 $$ 
 \mu = \frac{\rho}{\varepsilon} 
 \to 0. 
 $$ 
 
We seek the inner expansion in the form of the series 
\begin{equation}\label{ao} 
 H(\sigma,\omega,\mu)=\sum\limits_{n=0}^{\infty} 
 \mu^{n} h_n(\sigma,\omega), 
 \end{equation} 
 for whose coefficients from equation~(\ref{ea}) and
condition~(\ref{ic}) 
 we obtain the recurrence chain of initial value
problems 
\begin{align} 
 \frac{\partial h_0}{\partial \omega}-\frac{\partial^2
h_0}{\partial \sigma^2} &= 0, 
 & h_0(\sigma,0) &=\nu(\sigma), \label{ph0} \\ 
 \frac{\partial h_1}{\partial \omega}-\frac{\partial^2
h_1}{\partial \sigma^2} &= - \frac{\partial \varphi(h_0)}{\partial \sigma}, 
 & h_1(\sigma,0) &= 0, \label{ph1} \\ 
 \frac{\partial h_{n}}{\partial \omega}-\frac{\partial^2 h_{n}}{\partial \sigma^2} &= - \frac{\partial E_{n}}{\partial \sigma}, 
 & h_{n}(\sigma,0) &= 0, \label{phn} 
 \end{align} 
 where $$ 
 E_{n} = 
 \sum\limits_{q=1}^{n-1} 
 \frac{\varphi^{(q)}(h_0)}{q!} 
 \sum\limits_{n_1+\dotsc+n_q=n-1} 
 \prod\limits_{p=1}^{q} h_{n_p}, 
 \qquad 
 n\geqslant 2. 
 $$ 
 
To find where series~(\ref{ao}) makes sense
and to construct asymptotics in other domains,
 it is necessary to know the behavior
 of functions $h_n(\sigma,\omega)$ at infinity. 
 As shown in paper~\cite{hd},
the function 
\begin{equation}\label{fh0} 
 h_0(\sigma,\omega) = 
 \frac{1}{2\sqrt{\pi\omega}} 
 \int\limits_{-\infty}^{\infty} 
 \nu(s) \exp\left[ -\frac{(\sigma-s)^2}{4\omega}\right] 
 ds, 
 \end{equation} 
 i.e., the solution of problem~(\ref{ph0}),
  has the following asymptotics as $|\sigma|+\omega\to \infty$: 
 \begin{equation}\label{ah0} 
 h_0(\sigma,\omega) = R_{0,0,0}(z) + 
 \sum\limits_{m=1}^{\infty} \omega^{-m/2} 
 \left[ R_{0,m,0}(z) + \ln\omega R_{0,m,1}(z)\right], 
 \end{equation} 
 where 
\begin{equation}\label{r000} 
 R_{0,0,0}(z) = \nu^{-}_{0} 
\mathrm{erfc}(z) + \nu^{+}_{0} \mathrm{erfc}(-z), 
 \end{equation} 
 $$ 
 \mathrm{erfc} (z) = \frac{1}{\sqrt{\pi}} 
 \int\limits_{z}^{+\infty} \exp(-y^2)\,
dy, 
 \qquad 
 z = \frac{\sigma}{2\sqrt{\omega}}, 
 $$ 
 $R_{0,m,0}$ and $R_{0,m,1}$ are smooth
functions. 
 
 The asymptotics of solutions to problems~(\ref{ph1})--(\ref{phn}), 
 which can be expressed in the form of convolution 
\begin{equation}\label{hnr} 
 h_n(\sigma,\omega) = - \int\limits_{0}^{\omega} 
 \int\limits_{-\infty}^{\infty} 
 \frac{1}{2\sqrt{\pi(\omega-v)}} 
 \exp\left[ -\frac{(\sigma-s)^2}{4(\omega-v)}\right] 
 \frac{\partial E_n}{\partial s} 
 ds dv, 
 \end{equation} 
 as $|\sigma|+\omega\to \infty$ 
are found by the same method of paper~\cite{hd}. 
Proceeding by induction,
one can show that the following statement
is valid.

\textbf{Theorem~2.}
{\it
 For solutions of problems $(\ref{ph0})$--$(\ref{phn}),$ 
 which are determined recursively by formulas $(\ref{fh0})$ and $(\ref{hnr}),$ 
 for all $n\geqslant 0$
 there holds the asymptotic expansion 
\begin{equation}\label{ahn} 
 h_n(\sigma,\omega) = 
 \omega^{n/2} \sum\limits_{m=0}^{\infty} \omega^{-m/2} 
 \sum\limits_{l=0}^{m} (\ln\omega)^l 
 R_{n,m,l}(z). 
 \end{equation} 
}

\textbf{Acknowledgments.}
This work was supported by the Russian Foundation for Basic Research,
project no.~14-01-00322.

\end{document}